\documentclass[10pt]{amsart}
\usepackage{amsmath}
\usepackage{amssymb}
\usepackage{latexsym}

\newcommand{\Zint}{\mathbb {Z}}    
     
\newcommand{\R}{\mathbb {R}}      % Real number field
\newcommand{\halmos}{\rule{5pt}{5pt}}

\newtheorem{prop}{\bf Proposition}
\newtheorem{thm}[prop]{\bf Theorem}

\begin{document}
\title[Heun's differential equation and its $q$-deformation]
{Heun's differential equation and its $q$-deformation}
\author{Kouichi Takemura}
\address{Department of Mathematics, Faculty of Science and Engineering, Chuo University, 1-13-27 Kasuga, Bunkyo-ku Tokyo 112-8551, Japan}
\email{takemura@math.chuo-u.ac.jp}
\begin{abstract}
The $q$-Heun equation is a $q$-difference analogue of Heun's differential equation.
We review several solutions of Heun's differential equation and investigate polynomial-type solutions of $q$-Heun equation.
The limit $q\to 1$ corresponding to Heun's differential equation and the ultradiscrete limit $q\to 0$ are considered.
\end{abstract}
\maketitle

\section{Introduction}
Special functions have rich structures from the viewpoint of mathematics, and they have been applied to several branches of science.
A typical example of the special functions is the Gauss hypergeometric function and it is characterized as a solution of the Gauss hypergeometric differential equation, which is a standard form of the second order linear differential equation with three regular singularities $\{ 0,1, \infty \}$ on the Riemann sphere.

A standard form of second order linear differential equation with four regular singularities $\{ 0,1,t, \infty \}$ is given by
\begin{equation}
\frac{d^2y}{dz^2} +  \left( \frac{\gamma}{z}+\frac{\delta }{z-1}+\frac{\epsilon }{z-t}\right) \frac{dy}{dz} + \frac{\alpha \beta z -B}{z(z - 1)(z - t)} y= 0
\label{eq:Heun}
\end{equation}
with the condition $\gamma +\delta +\epsilon = \alpha +\beta +1$, and it is called Heun's differential equation.
The parameter $B$ is called an accessory parameter, which is independent from the local exponents.
Heun's differential equation and the differential equations of its confluent type appear in several systems of physics including quantum mechanics, general relativity, crystal transition and  fluid dynamics (e.g.~see \cite{Ron,SL,CH}).

A $q$-difference analogue of the Gauss hypergeometric function is called the basic hypergeometric function or $q$-hypergeometric function, and it has been studied vastly (see Gasper-Rahman \cite{GR}).
Hahn \cite{Hahn} introduced a $q$-difference analogue of Heun's differential equation of the form
\begin{equation}
%a(x) g(x/q) + b(x) g(x) + c(x) g(qx) =0
\{ a_2 x^2 +a_1 x+ a_0 \} g(x/q)  -\{ b_2 x^2 + b_1 x + b_0 \} g(x) + \{ c_2 x^2 + c_1 x+ c_0 \} g(xq) =0 ,
\nonumber
\end{equation}
where $a_2 a_0 c_2 c_0 \neq 0$. 
%such that $a(x)$, $b(x)$, $c(x)$ are polynomials such that $\deg_x a(x)= \deg_x c(x)=2 $, $a(0) \neq 0 \neq c(0) $ and $\deg _x b(x) \leq 2$.
It was rediscovered in \cite{TakR} by the fourth degeneration of Ruijsenaars-van Diejen system (\cite{vD0,RuiN}) or by specialization of the linear difference equation associated to the $q$-Painlev\'e VI equation (\cite{JS}).
Recently, the $q$-Heun equation appears in the study of the Heun-Askey-Wilson algebra (\cite{BTVZ}).

In this paper we review several solutions of Heun's differential equation and apply similar methods to some solutions of $q$-Heun equation.
Note that Heun's differential equation is recovered from the $q$-Heun equation by the limit $q\to 1 $ (see Eq.~(\ref{eq:Fuchs4sing}) and around that).

On the $q$-Heun equation, we may consider ultradiscrete limit $q\to +0$.
Then we expect thet the situation $q=0$ is simpler and it may help analysis of solutions of the $q$-Heun equation.

\section{A review of some solutions of Heun's differential equation.} \label{sec:revHeun}
First we investigate local solutions of Heun's differential equation about the regular singularity $z=0$.
Namely we investigate the solution to Eq.~(\ref{eq:Heun}) of the form $y= z ^{\mu } ( \tilde{c}_{0} + \tilde{c}_{1} z + \tilde{c}_{2} z^{2} + \cdots )$ where $\tilde{c}_{0}=1$.
By substituting it into Eq.~(\ref{eq:Heun}), we find that the value $\mu $ is either $0$ or $1-\gamma $.
These values are called exponents about the regular singularity $z=0$.
We can also calculate the exponents of other regular singularities and we can read them from the following Riemann scheme.
\begin{equation}
\left(
\begin{array}{cccc}
z=0 &  z=1 & z=t  & z=\infty \\
0 & 0 & 0 & \alpha \\ 
1- \gamma & 1- \delta & 1- \epsilon &  \beta 
\end{array}
\right)
\nonumber
\end{equation}
Then the parameter $B$ in Eq.~(\ref{eq:Heun}) is independent from the exponents and it is called the accessory parameter.

We now consider the solution of Eq.~(\ref{eq:Heun}) about $z=0$ of the exponent $0$, which is written as 
\begin{equation}
y=\sum_{n=0}^{\infty} \tilde{c}_{n} z^{n},\quad (\tilde{c}_{0}=1).
\label{eq:ysolHeun}
\end{equation}
By substituting it to Eq.~(\ref{eq:Heun}), the coefficients satisfy $t \gamma \tilde{c}_{1}   =B \tilde{c}_{0}$ and 
\begin{equation}
 t n(n-1+\gamma ) \tilde{c}_{n}  = [(n-1)\{(n-2+\gamma )(1+t)+  t \delta  + \epsilon \}+ B ] \tilde{c}_{n-1} -(n-2+\alpha )(n-2+\beta) \tilde{c}_{n-2}, \label{eq:Hlci}
% t (n+1)(n+\gamma ) \tilde{c}_{n+1}  = [n\{(n-1+\gamma )(1+t)+  t \delta  + \epsilon \}+ B ] \tilde{c}_{n} -(n-1+\alpha )(n-1+\beta) \tilde{c}_{n-1},  \label{eq:Hlci}
\end{equation}
for $n=2,3,\dots $.
If $t \neq 0,1$ and $\gamma \not \in {\mathbb{Z}}_{\leq 0}$, then the coefficients $\tilde{c}_{n}$ are determined recursively.
If we fix the accessory parameter $B$, the convergent radius of the formal solution $y$ in Eq.~(\ref{eq:ysolHeun}) is non-zero and no less than $\min (1, |t|)$.
The solution $y$ was denoted by $Hl(t,B;\alpha ,\beta ,\gamma ,\delta ;z)$ in \cite{Ron}.
Thus we can obtain a local solutions to Heun's differential equation.
Another solution can be obtained in the form $y= z ^{1- \gamma  } ( \tilde{c}'_{0} + \tilde{c}'_{1} z + \tilde{c}'_{2} z^{2} + \cdots ) $.

A main problem of Heun's differential equation is to investigate the value of the accessory parameter $B$ such that the differential equation admits a good global solution.

One type of a good solution is the non-zero solution which is simultaneously holomorphic about $z=0$ and $1$, which might be concerned with two-point boundary value problem.
For simplicity we assume that $|t|>1$.
Then the convergent radius of general local solution about $z=0$ is $1$.
For a special value of the accessory parameter $B$, the radius of convergence would be greater than $1$.
Then the solution is simultaneously holomorphic about $z=0$ and $1$.
The condition of the accessory parameter $B$ had been studied in several literature (e.g.~see \cite{Ron}), though the description is numerical.

Another type of a good solution is the polynomial solution.
We can check that a necessary condition for having a non-zero polynomial solution is that $\alpha$ or $\beta $ is a non-positive integer. 
Here we regard the accessory parameter $B$ to be an indeterminate.
Then the coefficient $\tilde{c}_{n}$ is a polynomial in $B$ of degree $n$ and we denote it by $\tilde{c}_{n}(B)$.

We assume that $\alpha=-N$ or $\beta=-N$ for some $N \in{\mathbb{Z}}_{\geq 0} $.
It follows from Eq.~(\ref{eq:Hlci}) for $n=N+2$ that 
\begin{equation}
 t (N+2) (N+1+\gamma ) \tilde{c}_{N+2} (B) = [(N+1)\{ (N +\gamma )(1+t) +t\delta +\epsilon \} +B ] \tilde{c}_{N+1} (B) + 0 \cdot \tilde{c}_N (B) .
\nonumber
\end{equation}
Let $B_{0}$ be a solution to the algebraic equation $\tilde{c}_{N+1} (B)=0$ of order $N+1$, i.e. $\tilde{c}_{N+1} (B_0)=0$.
Then we have $\tilde{c}_{N+2} (B_0)=0$ by the above equality.
By applying Eq.~(\ref{eq:Hlci}) for $n=N+2, N+3, \dots$, we have $\tilde{c}_{n} (B_{0})=0$ for $n\geq N+3$.
Hence, if $\tilde{c}_{N+1} (B_{0})=0$, then the differential equation (\ref{eq:Heun}) have a non-zero polynomial solution.
More precisely, we have the following proposition.
\begin{prop} $($\cite{WW,Ron}$)$ \label{prop:Heunpolym} 
Assume that $t \not \in \{ 0 ,1 \}$, $\gamma \not \in {\mathbb{Z}}_{\leq 0} $, $(\alpha +N)(\beta +N)=0$ and $N \in{\mathbb{Z}}_{\geq 0} $.
If $B$ is a solution to the equation $\tilde{c}_{N+1} (B)=0$, then the differential equation (\ref{eq:Heun}) have a non-zero polynomial solution of degree no more than $N$.
\end{prop}
We call $\tilde{c}_{N+1} (B)$ the spectral polynomial.
The polynomial solution $ y  = \tilde{c}_{0} +  \tilde{c}_{1} z + \cdots + \tilde{c}_{N} z^N $ is called the Heun polynomial.
Analysis of the spectral polynomial is one of the main problem of Heun's differential equation, and sometimes applied to physics.
Note that the Heun polynomials describe the eigenfunctions related with quasi-exact solvability.

Although the spectral polynomial $\tilde{c}_{N+1} (B)$ may have non-real roots or multiple roots in general, we have a sufficient condition that the spectral polynomial has only distinct real roots.
Namely, if $ (\alpha +N)(\beta +N)=0$, $N \in{\mathbb{Z}}_{\geq 0} $, $\delta $,$\epsilon $ and $\gamma $ are real, $\gamma >0$, $\delta + \epsilon  +\gamma +N >1$ and $t <0$, then the equation $\tilde{c}_{N+1} (B)=0$ has all its roots real and distinct.
Note that it can be proved by applying the argument of Sturm sequence (see \cite{CKLT2} for details).

\section{Solutions to $q$-Heun equation} \label{sec:qHeun}

We adopt the expression of the $q$-Heun equation as
\begin{eqnarray}
& (x-q^{h_1 +1/2} t_1) (x- q^{h_2 +1/2} t_2) g(x/q) + q^{\alpha _1 +\alpha _2} (x - q^{l_1-1/2}t_1 ) (x - q^{l_2 -1/2} t_2) g(qx) \label{eq:RuijD5} \\
& -\{ (q^{\alpha _1} +q^{\alpha _2} ) x^2 + E x + q^{(h_1 +h_2 + l_1 + l_2 +\alpha _1 +\alpha _2 )/2 } ( q^{\beta /2}+ q^{-\beta/2}) t_1 t_2 \} g(x) =0 \nonumber
\end{eqnarray}
(see \cite{TakR,TakqH}).
The parameter $E$ play the role of the accessory parameter.

We investigate the solution to Eq.~(\ref{eq:RuijD5}) of the form $y= x ^{\mu } ( c_{0} + c_{1} x + c_{2} x^{2} + \cdots )$ where $c_{0}=1$.
By substituting it into Eq.~(\ref{eq:RuijD5}), we find that the value $\mu $ is either $(h_1 +h_2 -l_1-l_2 -\alpha _1-\alpha _2 -\beta +2)/2$ or $(h_1 +h_2 -l_1-l_2 -\alpha _1-\alpha _2 +\beta +2)/2$.
Write
\begin{equation}
f(x)= x^{\lambda _1} \sum _{n=0}^{\infty } c_n x^n , \; \lambda _1 = \frac{h_1 +h_2 -l_1-l_2 -\alpha _1-\alpha _2 -\beta +2}{2},
\label{eq:A4la1la2}
\end{equation}
and substitute it into Eq.~(\ref{eq:RuijD5}).
Then the coefficients $c_{n} $ $(n=1,2,\dots )$ are determined recursive by 
\begin{eqnarray}
& & c_n t_1 t_2 [ q^{h_1 +h_2 } ( 1 - q^{n }) ( 1 - q^{n -\beta }) ] \label{eq:rec01}  \\
& & = c_{n-1} [E q^{n -1 +\lambda _1} + q^{1/2}(q^{h_1 } t_1 +q^{h_2 } t_2 ) +  (q^{l_1 } t_1 +q^{l_2 } t_2 ) q^{2(n +\lambda _1) +\alpha _1 +\alpha _2 -5/2 }]  \nonumber \\
& & \qquad - c_{n-2} [ q (1 - q^{n -2 +\lambda _1 +\alpha _1})(1 - q^{n-2 +\lambda _1 + \alpha _2}) ] , \nonumber 
\end{eqnarray}
with the initial condition $c_0=1$ and $c_{-1} =0$ (see \cite{TakqH,KST}).
If we regard $E$ as an indeterminant, then $c_n $ is a polynomial of $E $ of degree $n$, and we denote it by $c_n (E)$

If $\lambda _1 +\alpha _1$ is a non-positive integer, then $1 - q^{n -2 +\lambda _1 +\alpha _1} =0 $ for $n= -\lambda _1 -\alpha _1 +2$ and we may apply a similar argument of the polynomial solutions of Heun's differential equation.
The polynomial-type solution of the $q$-Heun equation, which is written as a terminating series, is described as follows.
\begin{thm} $($\cite{KST}$)$ \label{prop:prop0}
Let $\lambda _1$ be the value in Eq.~(\ref{eq:A4la1la2}) and assume that $-\lambda _1 - \alpha _1 (:=N)$ is a non-negative integer and $\beta \not \in \{ 1,2,\dots ,N, N+1\}$.
Set $c_{-1}(E)=0 $,  $c_0(E)=1 $ and we determine the polynomials $c_n(E)$ $(n=1,\dots ,N+1)$ recursively by Eq.~(\ref{eq:rec01}). 
If $E=E_0$ is a solution of the algebraic equation
\begin{equation}
c_{N+1} (E)=0 ,
\nonumber
%\label{eq:cN+1}
\end{equation}
then $q$-Heun equation defined in Eq.~(\ref{eq:RuijD5}) has a non-zero solution of the form
\begin{equation}
f(x)= x^{\lambda _1} \sum _{n=0}^{N } c_n (E_0)  x^n .
\label{eq:gxxlapol}
\end{equation}
\end{thm}
We call $c_{N+1} (E)$ the spectral polynomial of the $q$-Heun equation, and we call $f(x)$ in Eq.~(\ref{eq:gxxlapol}) the polynomial-type solution, which is a product of $x^{\lambda _1} $ and a polynomial.

We may use the theory of Sturm sequence from the three term relations in Eq.~(\ref{eq:rec01}) and we obtain real root property of the spectral polynomial $c_{N+1} (E) $.
\begin{thm} $($\cite{KST}$)$ \label{thm:Sturm}
Assume that $N= -\lambda _1 -\alpha _1 $ is a non-negative integer, $t_1$, $t_2$, $h_1$, $h_2$, $l_1$, $l_2$, $\alpha _1$, $\alpha _2$, $\beta  $ are all real, $t_1 t_2 >0 $ and $q$ is a positive number such that $q \neq 1$.
If $\alpha _2-\alpha _1<1 $ and $ \beta >-1$, then the equation $c_{N+1} (E)=0$ has all its roots real and distinct.
\end{thm}
By the theorem, we can label the roots of the spectral polynomial $c_{N+1} (E) $ as $E_1 (q), E_2(q), \dots ,E_{N+1} (q)$ such that $E_1 (q) < E_2(q) < \dots < E_{N+1} (q) $ for $0<q<1$.
Note that the real root property of the spectral polynomial for several other cases were shown in \cite{KST}.

\section{Limit to Heun's differential equation} \label{sec:limq1}

We consider the limit $q \to 1$ for the $q$-difference equation and its solutions.
For the $q$-difference equation, the limit $q \to 1$ was discussed in \cite{TakR} and we review it here.
Set $q= 1+ \varepsilon $ and replace the accessory parameter $E$ to $B$ by the relation
\begin{equation}
E= - 2 (t_1+t_2) -\varepsilon C_1 +\varepsilon ^2 (B - C_2) , \nonumber
\end{equation}
where
\begin{eqnarray}
& & C_1= (\alpha _1 +\alpha _2 )(t_1+t_2) +(l_1 +h_1 )t_1 +(l_2 +h_2 )t_2, \nonumber \\
& & C_2= \frac{t_1}{2} \Big\{ h_1^2 + (\alpha _1 +\alpha _2+ l_1 - 1)^2 - \frac{1}{2} \Big\} + \frac{t_2}{2} \Big\{ h_2^2 + (\alpha _1 +\alpha _2+ l_2 - 1)^2 -\frac{1}{2} \Big\}.  \nonumber 
\end{eqnarray}
By using Taylor's expansion
\begin{eqnarray}
& & g(x/q) =g(x) + (-\varepsilon +\varepsilon ^2 )xg'(x) + \varepsilon ^2 x^2 g''(x) /2 +O( \varepsilon ^3), \nonumber \\
& & g(qx) =g(x) + \varepsilon  x g'(x) + \varepsilon ^2 x^2 g''(x) /2 +O( \varepsilon ^3), \nonumber 
\end{eqnarray}
we rewrite Eq.~(\ref{eq:RuijD5}) as a series of $\varepsilon   $.
Then the coefficients of $\varepsilon ^0$ and $\varepsilon ^1 $ are cancelled and the coefficient of $\varepsilon ^2 $ is written as 
\begin{eqnarray}
& & x^2(x-t_1)(x-t_2) g''(x) \label{eq:Fuchs4sing} \\
& & +  [ (1+ h_2 - l_2)x(x- t_1 ) + (1 +h_1 - l_1 )x(x- t_2) + ( \tilde{l}  -1 ) (x -t_1)(x-t_2) ]  xg'(x) \nonumber \\
& & + [ \alpha _1 \alpha _2 x^2 - B x + t_1 t_2 ( \tilde{l}/2 -1 +\beta/2)( \tilde{l}/2 -1 -\beta/2)] g(x) =0, \nonumber 
\end{eqnarray}
where $\tilde{l} = l_1 +l_2 +\alpha _1 +\alpha _2 -h_1 - h_2  $ (see also \cite{TakR}).
Heun's differential equation was obtained from Eq.~(\ref{eq:Fuchs4sing}) by setting $z=x/t_1$ and $g(x)=x^{1-\tilde{l}/2 - \beta/2 }\tilde{g}(z)$ in \cite{TakR}.
Here we obtain Heun's differential equation by restricting the parameters. 
Set $ \tilde{l}/2 -1 +\beta/2 =0$ $(\Leftrightarrow \lambda _1=0)$, $t_1=1$ and $t_2=t$. 
Then we obtain Heun's differential equation
\begin{equation}
\frac{d^2y}{dx^2} +  \left( \frac{\gamma}{x}+\frac{\delta }{x-1}+\frac{\epsilon}{x-t}\right) \frac{dy}{dx} + \frac{\alpha _1 \alpha _2 x -B}{x(x - 1)(x - t)} y= 0,
\nonumber
\end{equation}
where $\gamma = 1 - \beta$, $\delta = 1 +h_1 -l_1$ and $\epsilon = 1+h_2-l_2$.

We consider the limit of the coefficients of the local expansion in Eq.~(\ref{eq:A4la1la2}).
By the limit $ \varepsilon \to 0$ $(\Leftrightarrow q \to 1)$, the recursive relation Eq.~(\ref{eq:rec01}) tends to the following equation
\begin{eqnarray}
& & c_n n (n-\beta ) t_1 t_2 = c_{n-1} [ B + (n+ \lambda _1 -1) \{ (n+\lambda _1 -h_1+l_1+\alpha _1+\alpha _2-2) t_1 \label{eq:recvep0} \\
& & \qquad + (n+\lambda _1 -h_2+l_2+\alpha _1+\alpha _2-2) t_2 \} ]  - c_{n-2} [  (n -2 +\lambda _1 +\alpha _1)(n-2 +\lambda _1 + \alpha _2 ) ] . \nonumber 
\end{eqnarray}
If $\beta $ is not positive integer, then the coefficients $c_n $ $(n=1,2,\dots  )$ are determined recursively by setting $c_0=1 $ and $c_{-1} =0$, and $c_n$ is a polynomial of $B$ of degree $n$.
We denote $c_n$ by $\tilde{c}_n (B)$.
To discuss polynomial-type solutions, we assume that $-\lambda _1 - \alpha _1$ is a non-negative integer and $\beta \not \in \{ 1,2,\dots , N+1\}$, where $N= -\lambda _1 - \alpha _1$.
As Theorem \ref{prop:prop0}, we can show that Eq.~(\ref{eq:Fuchs4sing}) has a non-zero solution of the form $f(x)= x^{\lambda _1} \{ 1 + c_1 x + \cdots + c_N x^N \} $, if $B=B_0$ is a solution of the algebraic equation $\tilde{c}_{N+1} (B)=0 $.
We can obtain real root property of the spectral polynomial $\tilde{c}_{N+1} (B) $ by using the theory of Sturm sequence.
Namely, if $t_1$, $t_2$, $h_1$, $h_2$, $l_1$, $l_2$, $\alpha _1$, $\alpha _2$, $\beta  $ are all real, $\alpha _2-\alpha _1<1 $ and $ \beta >-1$, then the equation $\tilde{c}_{N+1} (B)=0$ has all its roots real and distinct.
Hence we can label the roots of the spectral polynomial $\tilde{c}_{N+1} (B) $ as $B_1 , B_2, \dots ,B_{N+1} $ such that $B_1 < B_2 < \dots < B_{N+1} $.
Combining with the limit procedure, we have 
\begin{equation}
B_j= C_2 + \lim _{q \to 1-0} \frac{E _j (q) + 2 (t_1+t_2)  + (q-1) C_1 }{(q-1)^{2}}.
\nonumber
\end{equation}

We restrict the parameters in Eq.~(\ref{eq:recvep0}) to discuss with Heun's differential equation.
Set $\lambda _1 =0 $, $t_1=1$ and $t_2=t$.
Then it follows from $\lambda _1 =0 $ that $\alpha _1 +\alpha _2 = h_1 +h_2 -l_1-l_2  -\beta +2 $ and
\begin{eqnarray}
& n (n-\beta ) t c_n  =[ B + (n -1) \{ (n +h_2 -l_2 -\beta )
+ (n +h_1 -l_1 -\beta ) t \} ]  c_{n-1} \nonumber \\
& \qquad  - (n -2 +\alpha _1)(n-2 + \alpha _2 ) c_{n-2}  . \nonumber 
\end{eqnarray}
Hence we recover Eq.~(\ref{eq:Hlci}) by setting $\gamma = 1 - \beta$, $\delta = 1 +h_1 -l_1$, $\epsilon = 1+h_2-l_2$ and replacing $\alpha $ and $\beta $ with $\alpha _1 $ and $\alpha _2 $.

\section{Analysis of the spectral polynomial by the ultradiscrete limit} \label{sec:UDL}

In general it would be impossible to solve the roots of the spectral polynomial $c_{N+1} (E)$ of $q$-Heun equation and those of Heun's differential equation explicitly.
To understand the roots of the spectral polynomial, we may apply the ultradiscrete limit $q \to +0$.
In \cite{KST}, the behaviour of the roots $E=E_1 (q), E_2 (q) , \dots ,E_{N+1} (q) $ of $c_{N+1} (E)=0 $ by the ultradiscrete limit was studied in some cases.

As discussed in \cite{KST}, we define the equivalence $\sim $ of functions of the variable $q$ by 
\begin{equation}
a(q) \sim b (q) \; \Leftrightarrow \; \lim _{q \to +0} \frac{a(q)}{b(q)} =1.
\nonumber
\end{equation}
We also define the equivalence $\sum _{j=0}^{M} a_j(q) E^j \sim \sum _{j=0}^{M} b_j (q) E^j $ of the polynomials of the variable $E$ by $ a_j (q) \sim b _j(q) $ for $j=0,\dots ,M$.
If $\sum _{j=0}^{M} a_j(q) E^j \sim \sum _{j=0}^{M} c_j q^{\mu _j} E^j $ for some $c_j \in \R \setminus \{ 0\}$ $(j=0,1,\dots ,M)$, then we call $\sum _{j=0}^{M} c_j q^{\mu _j} E^j $ the leading term of $\sum _{j=0}^{M} a_j(q) E^j $.

We are going to find simpler forms of the polynomials $c_{n} (E)$ $(n=1,2,\dots ,N+1)$ determined by Eq.~(\ref{eq:rec01}) with respect to the equivalence $\sim $. 
For simplicity, we assume 
\begin{equation}
0<q<1, \; N= -\lambda _1 -\alpha _1 \in \Zint _{\geq 0}, \; \beta <1, \; \alpha _2-\alpha _1<1, \; t_1>0 , \; t_2>0, \; h_1 <h_2 , \; l_1<l_2,
\label{eq:assump}
\end{equation} 
which was also assumed in \cite{KST,KST2}.
The condition in Eq.~(\ref{eq:assump}) implies the condition in Theorem \ref{thm:Sturm} such that the roots of the spectral polynomial $c_{N+1} (E) $ are real and distinct. 
It was shown in \cite{KST} that if $1 +h_2 -l_2 - \beta >0 $ and the condition in Eq.~(\ref{eq:assump}) is satisfied, then the three term relation in Eq.~(\ref{eq:rec01}) is simplified as
\begin{equation}
 c_n (E) \sim t_1^{-1} t_2^{-1} [ E q^{n -1 -h_1 -h_2 +\lambda _1} + t_1 q^{1/2 -h_2 } ] c_{n-1} (E) - t_1^{-1} t_2^{-1} q^{2n-1 -l_1-l_2 -\beta } c_{n-2} (E)
\nonumber
\end{equation}
for $ n=1,2,\dots ,N+1$ under the assumption that there are no cancellation of the leading terms of the coefficients of $E^j$ $(j=0,1,\dots ,n-1)$ in the right hand side.
By adding one more condition, we may ignore the term $t_1^{-1} t_2^{-1} q^{2n-1 -l_1-l_2 -\beta } c_{n-2} (E) $ and we have the following theorem.
\begin{thm} $($\cite{KST}$)$ \label{thm:specpol1}
We assume Eq.~(\ref{eq:assump}), $1 +h_2 -l_2 - \beta >0 $ and $2+ 2 h_2  - l_1 - l_2 -\beta >0  $.\\
(i) The spectral polynomial $c_{N+1} (E)$ satisfies
\begin{eqnarray}
& c_{N+1} (E) \sim  (t_1 t_2 )^{-N-1} q^{(N/2 + \lambda _1 -h_1 -h_2)(N+1) } ( E + q^{1/2 -N +h_1 -\lambda _1} t_1 )( E + q^{3/2 -N +h_1 -\lambda _1} t_1 ) \nonumber \\
& \qquad \cdots ( E + q^{-1/2 +h_1 -\lambda _1 } t_1 )( E + q^{1/2 +h_1 -\lambda _1 } t_1 ) . \nonumber 
\end{eqnarray}
(ii) The roots $E _k(q)$ $(k=1, \dots ,N+1)$ of the spectral polynomial $c_{N+1} (E) $ such that $E _1(q) < E _2(q)< \dots < E _{N+1} (q) $ satisfies
\begin{eqnarray}
& E _k(q)  \sim -q^{k -N -1/2 +h_1 -\lambda _1  } t_1 .
\nonumber
\end{eqnarray}
\end{thm}
We can investigate the asymptotic of the zeros of the Heun polynomial in Eq.~(\ref{eq:gxxlapol}) as $q \to +0$ on the situation of Theorem \ref{thm:specpol1}.
\begin{thm} $($\cite{KST}$)$ \label{thm:eigenf1}
Let $k\in \{ 1,2,\dots ,N+1 \}$.
Assume Eq.~(\ref{eq:assump}), $1 +h_2 -l_2 - \beta >0 $, $2+ 2 h_2  - l_1 - l_2 -\beta >0  $ and the value $E=E_{N+2-k} (q) $ is a solution of the characteristic equation $c_{N+1} (E)=0$ such that $E_{N+2-k} (q)  \sim -q^{3/2 -k +h_1 -\lambda _1} t_1 $.
The Heun polynomial $\sum _{n=0}^{N } c_n (E_{N+2-k} (q))  x^n  $ has zeros $x= x_j(q) $ $(j=1,\dots ,N)$ for sufficiently small $q$, which are continuous on $q$ and satisfy
\begin{equation}
 \lim _{q \to +0} \frac{x _j (q)}{q^{j -1/2 +h_2 }t_2 } =1 , \; (j=1,\dots ,k-1), \quad  \lim _{q \to +0} \frac{x _j (q)}{- q^{-2j -1/2 -h_2 +l_1 +l_2 +\beta }t_1 } =1 ,\; (j=k , \dots ,N). 
\nonumber
\end{equation}
\end{thm}
Note that several other cases of the ultradiscrete limit of the spectral polynomial were studied in \cite{KST,KST2}.

\section*{Acknowledgements}
The author was supported by JSPS KAKENHI Grant Number 18K03378.

\end{document}